\documentclass[preprint,12pt,3p]{elsarticle}

\usepackage{amssymb}

\usepackage{amsmath}
\journal{Arxiv}

\begin{document}

\begin{frontmatter}

\title{ON THE GENERATION OF  PYTHAGOREAN TRIPLES AND REPREZENTATION OF INTEGERS AS A DIFFERENCE OF TWO SQUARES}

\author{Emil Asmaryan\corref{cor1}\fnref{label2}}
\address[label2]{Institute of  Radiophysics and Electronics of Armenian Nat.Ac.Sci., Alikhanian Bros.str.,1, 0203, Ashtarak, Armenia}

\cortext[cor1]{Emil Asmaryan}

\ead{emilas@irphe.am}



\begin{abstract}
The general formulas for finding the quantity of all primitive and nonprimitive triples generated by the given number $x$ have been proposed. Also the formulas for finding the complete quantity of the representations of the integers  as a difference of  two squares  have  been  obtained.
\end{abstract}

\end{frontmatter}


\section{Introduction}\label{sec1}
The equation 
\begin{equation}\label{eq:1}
x^2+y^2=z^2	
\end{equation}
is under consideration. The solutions of \eqref{eq:1}, when $x$,$y$,$z$ ($x$,$y>2$) are the integer positive numbers, are named the Pythagorean triples, and \eqref{eq:1} itself ranks among the Diofantine equation. If some triple ($x$,$y$,$z$) is known, then one can obtain the infinite quantity of solutions of \eqref{eq:1} by multiplying of  $x$,$y$,$z$  by the arbitrary numbers. The triples with coprime $x$,$y$,$z$ are named primitive.

Since the  Euclidean times there is a well-known method, which allows to calculate the Pythagorean triples by two integer numbers $k$, $l$ , representing $x$, $y$, $z$ as

\begin{equation*}
x=2kl, \hspace{1cm} y=k^2-l^2, \hspace{1cm} z=k^2+l^2; \hspace{1cm} k>1 
\end{equation*}
\begin{equation}\label{eq:2}
(2kl)^2+(k^2-l^2)^2=( k^2 + l^2)^2
\end{equation}
If  $k$ and  $l$ are coprime and of opposite parity, then \eqref{eq:2} gives the primitive triples.

There also exist the other ways of obtaining of the Pythagorean triples or their representations as special combinations of certain numbers.

However, in this work the problem of finding the Pythagorean triples is stated as follows:
for the arbitrary integer $x>2$ to find the total quantity of the primitive and nonprimitive triples, in which $x$ is one of the elements on the left side of \eqref{eq:1} (i.e. generated by $x$), and also to propose the convenient way for their calculation.

The methods mentioned above are not effective to solve this kind of problem. So, the method \eqref{eq:2} gives only as many triples as the given $x$ can be represented as $2kl$ (if $x$ is even), or as $k^2-l^2$  (if  $x$ is odd). 

The Euclidean method gives the right quantity of primitive triples  and also gives some quantity of nonprimitive (if the integer numbers  $k$ and $l$ are not coprime, or they are both odd). However, as it is shown below in this paper, the total quantity  of nonprimitive triples, generated by $x$, in reality may be significantly more than the quantity of the triples obtained from \eqref{eq:2} with the integers $k$ and $l$. 

In the present work the general formulas for finding the quantity of all triples generated by the given number $x$ have been proposed. Also the formulas for finding the complete quantity of the representations of the integers  as a difference of  two squares have  been  obtained. Some of  this results are represented in \cite{proc2016}.

\section{Calculation of triples}

Let us turn to the calculation of  Pythagorean triples and their quantity.
Taking into account, that the numbers $x$ and  $y$  can both be even or of different parity, but cannot  both be odd, we use the following standard procedure.
\begin{enumerate}[(a)] 
	\item \label{en:a} Let the odd number $x >1$ be  given. Then $y$ is  even, and $z$ is odd, i.e. $z=y+(2m+1)$, $m \geq 0$.
	 
	From \eqref{eq:1} we have:
		
	\begin{equation*}
     x^2+[z-(2m+1)]^2=z^2,
     \end{equation*}
     \begin{equation}\label{eq:3}
     z=\frac{1}{2}\bigg[\frac{x^2}{2m+1}+2m+1\bigg],\hspace{1cm} y=\frac{1}{2}\bigg[\frac{x^2}{2m+1}-(2m+1)\bigg].
	\end{equation}
	
	Since $y$ and $z$ must be integer, then $2m +l\equiv d$ must be divisor of $x^2$:
	\begin{equation}\label{eq:4}
	y_i=\frac{1}{2}\bigg(\frac{x^2}{d_i}-d_i\bigg),\hspace{1cm} z_i=\frac{1}{2}\bigg(\frac{x^2}{d_i}+d_i\bigg)
	\end{equation}
	Substituting into \eqref{eq:4} the values of divisors $d_i$, we obtain corresponding triples ($x$, $y_i$, $z_i$).
	
	It is easy to show, that from \eqref{eq:4} just $z$ will be odd, and $y$ will be even.
	
	In order to find the quantity of all different triples with positive elements, only divisors $d_i<x$  must be used in \eqref{eq:4}. Let  $N_d$ be the quantity of all divisors of $x^2$. It is obvious that the  quantities of divisors  $d_i<x$ and $d_i>x$ are  the same. Therefore the quantity of $d_i<x$ and so the quantity of Pythagorean triples is:
	\begin{equation}\label{eq:5}
	N_{tr}=\frac{N_d-1}{2}
	\end{equation}
	Note, that the squares  of integer numbers always have the odd quantity of  divisors. 
	\item \label{en:b}	Let the even number $x>2$ be given. In this case $y$ and $z$ are the numbers of the same parity, i.e. $z=y+2m$, $m>0$.
	
	Then, from \eqref{eq:1} we obtain:
	\begin{eqnarray*}
	x^2+(z-2m)^2=z^2,\\
	x^2-4mz+4m^2=0
	\end{eqnarray*}
	\begin{equation}\label{eq:6}
	y=\frac{(\frac{x}{2})^2}{m}-m,\hspace{1cm} z=\frac{(\frac{x}{2})^2}{m}+m
	\end{equation}
	Since $y$ and $z$ must be integer, then $m\equiv d$ must be divisor of  $(\frac{x}{2})^2$ ,i.e.
	\begin{equation}\label{eq:7}
	y_i=\frac{(\frac{x}{2})^2}{d_i}-d_i,\hspace{1cm} z_i=\frac{(\frac{x}{2})^2}{d_i}+d_i.
	\end{equation}
	Substituting into\eqref{eq:7}  the values of $d_i$, we obtain the triples ($x$, $y_i$, $z_i$). In order to obtain the triples with positive elements only divisors $d_i<\frac{x}{2}$   must be used in \eqref{eq:7}. If  $N_d$ is the quantity  of  all  divisors of $(\frac{x}{2})^2$, then the quantity of divisors $d_i<\frac{x}{2}$ is equal to  $\frac{N_d-1}{2}$. Therefore, the quantity of all  Pythagorean triples, generated by the number $x$, is also $\frac{N_d-1}{2}$.

	Let us  now determine  $N_d$ in  cases \ref{en:a}) and  \ref{en:b}). For convenience we redenote $x\equiv n$.
	
	As it is known, the complete quantity $Q$ of divisors of any number $n$ may be determined by its canonic expansion
	 
	\begin{equation}\label{eq:8}
	n=p_1^{s_1}\cdot p_2^{s_2}\cdot\cdot\cdot \cdot p_q^{s_q},
	\end{equation}
	where  $p_1,\cdot\cdot\cdot,p_q$ are the different prime numbers none-equal to 1.
	
\end{enumerate}

Then 
\begin{equation}\label{eq:9}
Q=(s_1+1)(s_2+1)...(s_q+1).
\end{equation} 

If  $n$  is odd number,  then from  its canonic expansion \eqref{eq:8} we obtain:

\begin{equation}\label{eq:10}
n^2=p_1^{2s_1}\cdot p_2^{2s_2}\cdot\cdot\cdot \cdot p_q^{2s_q},
\end{equation}
and
\begin{equation}\label{eq:11}
N_d=(2s_1+1)(2s_2+1)...(2s_q+1).
\end{equation}
Hence the complete quantity of  Pythagorean triples is:
\begin{equation}\label{eq:12}
N_{tr}=\frac{(2s_1+1)(2s_2+1)...(2s_q+1)-1}{2}
\end{equation}
If  $n$  is even, then $p_1=2$.  Writing down $n$ as  $n=2^{s_1+1}\cdot p_2^{s_2}\cdot\cdot\cdot p_q^{s_q}$, $s_1\geq 0$,  we obtain: 
\begin{equation}\label{eq:13}
\frac{n}{2}=2^{s_1}\cdot p_2^{s_2}\cdot\cdot\cdot p_q^{s_q},\hspace{1cm}
(\frac{n}{2})^2=2^{2s_1}\cdot p_2^{2s_2}\cdot\cdot\cdot p_q^{2s_q}
\end{equation}
Then  the complete quantity of triples is given by the same expression \eqref{eq:12}, but $s_1,s_2,\cdot\cdot\cdot,s_q$ are taken from the canonic expansion \eqref{eq:13} of   $\frac{n}{2}$:
Thus we can formulate the following result: 

Every integer number $n>2$  generates 〖$N_{tr}$ Pythagorean triples  having the form:
\begin{enumerate}[(a)]
	\item \label{en:a2}	for odd n:
	\begin{equation}\label{eq:14}
	\left\{ n,\frac{1}{2}\bigg(\frac{n^2}{d_i}-d_i\bigg),\frac{1}{2}\bigg(\frac{n^2}{d_i}+d_i\bigg) \right\},
	\end{equation}
	where $d_i$  are the divisors of $n^2$ which are less than  $n$,  and  〖$N_{tr}$ is determined by expression \eqref{eq:12} with $s_1,s_2,\cdot\cdot\cdot,s_q$  taken from canonic expansion  of  $n$;
	\item \label{en:b2}	for even  $n$:
	\begin{equation}\label{eq:15}
	\left \{ n,\frac{(\frac{n}{2})^2}{d_i}-d_i,\frac{(\frac{n}{2})^2}{d_i}+d_i\right\},
	\end{equation}
	where $d_i$  are the divisors of $(\frac{n}{2})^2$ which are less than  $\frac{n}{2}$,  and  〖$N_{tr}$ is determined by expression  \eqref{eq:12} with   $s_1,s_2,\cdot\cdot\cdot,s_q$  taken from canonic expansion  of $\frac{n}{2}$.
	
	As the illustration let us consider the example:
$n=120= 2^3\cdot3\cdot5;   \frac{n}{2}= 2^2\cdot3\cdot5=60;   (\frac{n}{2})^2= 3600 = 2^4\cdot3^2\cdot5^2;   N_d= 45; N_{tr}=22,$ including 4 primitive. Euclidean method \eqref{eq:2} (with integer $k$, $l$) gives only 6 triples, including 4 primitive.

\end{enumerate}

It follows from \eqref{eq:14} and \eqref{eq:15} that for given $n$ there are maximal and minimal values of other two elements of Pythagorean triples, namely:

\begin{enumerate}[(a)]
	\item \label{en:a3}	for odd n:\newline
	maximal:$ \frac{n^2-1}{2}$    and $\frac{n^2+1}{2}$, at  $d_i = 1$ \newline
	minimal: $\frac{1}{2}(\frac{n^2}{d_{max}}-d_{max})$ and $\frac{1}{2}(\frac{n^2}{d_{max}}+d_{max}),$ where $d_{max}$ is the divisor of $n^2$, most near to $n$;
	
	\item \label{en:b3}	for even n:\newline
	maximal:$ \big(\frac{n}{2}\big)^2-1$    and $\big(\frac{n}{2}\big)^2+1$, at  $d_i = 1$ \newline
	minimal: $\frac{\big(\frac{n}{2}\big)^2}{d_{max}}-d_{max}$ and $\frac{\big(\frac{n}{2}\big)^2}{d_{max}}+d_{max}$, where  
	$d_{max}$ is the divisor of $(\frac{n}{2})^2$, most near to  $\frac{n}{2}$.	
\end{enumerate}

\section{Determination of the quantity of primitive triples}
Now we find the quantity of primitive triples included in $N_{tr}$  
\subsection{Let number $n$ be odd}
We use again the canonic expansion \eqref{eq:10}. As a divisors  $d_i$ in \eqref{eq:14} let’s take multipliers contained in \eqref{eq:10}, which are less than $n$ and represent the products of numbers $p_j^{2s_j}$ between themselves in all possible quantities and combinations (single, double, triple etc), i.e. the numbers of form 
\begin{eqnarray}\label{eq:16}
p_j^{2s_j},p_j^{2s_j}\cdot p_m^{2s_m},p_j^{2s_j}\cdot p_m^{2s_m}\cdot p_j^{2s_r},... < n \\
j,m,r...=1,2,...,q \nonumber
\end{eqnarray}

We see, that when dividing  the number $n^2$ , written in form \eqref{eq:10}, by the divisors $d_i$ of the type \eqref{eq:16}, all appropriate $p_j,...,p_r$ will be cancelled in the first term of bracket in expressions $\frac{1}{2}\big(\frac{n^2}{d_i}\mp d_i\big)$, whereas the second term contains only these $p_j,...,p_r$. Therefore this two terms will not contain the same $p_j,...,p_r$  and cannot divide by any divisors of $n$. Then the elements of such triples will be coprime and these triples will be  primitive.

Since the numbers $p_j^{2s_j}$ are contained  in $d_i$ giving the primitive triples in indivisible manner (i.e. they can be considered as if  they are “prime”, in power 1), then writing \eqref{eq:10} as
\begin{equation*}
n^2=\big( p_1^{2s_1}\big)^1\cdot....\cdot \big(p_q^{2s_q}\big)^1
\end{equation*}

and  using \eqref{eq:9}, we obtain that the quantity of divisors of type \eqref{eq:16}is equal to $\frac{(1+1)^q}{2}=2^{q-1}$. Hence, we obtain the following result:

The quantity of primitive triples generated by the odd number $n$, is: 
\begin{equation}\label{eq:17}
N_p=2^{q-1}
\end{equation}

where $q$ is the quantity of prime numbers $p_j$ in the canonic expansion of $n$.

Thus, $N_p$ does not depend on numbers  $p_j,s_j$  and is determined only by number $q$.
\subsection{Let $n$ be even} If $\frac{n}{2}$  is odd (i.e. $n$ indivisible by 4), then one can see from \eqref{eq:15} that $y$ and $z$ are also both even by all divisors $d_i$. Therefore such $n$ cannot generate the primitive triples. If $\frac{n}{2}$ is even, then the primitive triples are obtained from \eqref{eq:15}  by $d_i$, having the form \eqref{eq:16}  in canonic expansion \eqref{eq:13} of $\big(\frac{n}{2}\big)^2$. Hence, we obtain the following combined expression for quantity of primitive triples generated by even $n$:
\begin{equation}\label{eq:18}
N_p=2^{q-2}\big[1+(-1)^{n/2}\big]
\end{equation}
where $q$ is the quantity of prime numbers $p_j$ in canonic expansion of  $\frac{n}{2}$.
 
Note that $N_p$  always includes the triple, obtained by $d_i =1$.

For illustration we consider the example:
\begin{equation*}
n=2220=2^2\cdot 3 \cdot 5\cdot 37; \frac{n}{2}=1110=2\cdot 3\cdot 5 \cdot 37;  \bigg(\frac{n}{2}\bigg)^2=2^2\cdot 3^2 \cdot 5^2 \cdot 37^2 =1232100
\end{equation*}

The total quantity of triples  $N_{tr} = \frac{(2+1)^4-1}{2}=40$;

The quantity of primitive triples $N_p=2^{4-1} = 8$, and they are obtained from \eqref{eq:15} by the following divisors $d_i$:
$d_1 = 1; d_2= 2^2; d_3 = 3^2; d_4 = 5^2; d_5 = 2^2\cdot 3^2; d_6 = 2^2\cdot 5^2; d_7 = 3^2\cdot 5^2; d_8 = 2^2\cdot 3^2\cdot 5^2$.

In works \cite{robbins2006number,OMLAND20171} the quantity of primitive Pythagorean triangles with given inradius has been obtained. Here we have determined the complete quantity of Pythagorean triangles with given cathetus.

\section{The representation of integers as a difference of the  squares of two integer numbers}

The formulas \eqref{eq:14} and \eqref{eq:15} are proven to be useful for finding all possible representations of integer  numbers  as a difference of two squares.
\subsection{}  
Let $n$ be the arbitrary odd number more than 1, and we want to represent it in form
\begin{equation}\label{eq:19}
n=k^2-l^2, \hspace{1cm} k,l>0
\end{equation}
It is clear from identity \eqref{eq:2} written as 
\begin{equation*}
	(2kl)^2+n^2=(k^2+l^2)^2,
\end{equation*}
that $2kl$  and  $k^2+l^2$ are the elements of  Pythagorean triples generated by odd number $n$. As we already know, they are given by expressions \eqref{eq:14}.

In particular, we can write:

\begin{equation}\label{eq:20}
k^2+l^2=\frac{1}{2}\bigg(\frac{n^2}{d_i}+d_i\bigg)
\end{equation}

From \eqref{eq:19} and \eqref{eq:20} we find $k$ and $l$:
\begin{equation}\label{eq:21}
k=\frac{n+d_i}{2\sqrt{d_i}}, \hspace{1cm} l=\frac{n-d_i}{2\sqrt{d_i}}
\end{equation}

It follows from \eqref{eq:21}, that $k$  and $l$ will be integers and positive, if  $d_i< n$ and are the squares of integer numbers. Indeed, in this case $\sqrt{d_i}$   is the integer number, being the divisor of  $n$  and  $d_i$,  and, besides, since $n$  and  $d_i$,  are both odd, so their sum and difference are even. One can see that $k$  and $l$  have the opposite parity.

Therefore, the quantity $N_r$  of all representations \eqref{eq:19} is equal to quantity of divisors of $n^2$, which are less than $n$, and are in canonic expansion \eqref{eq:10} the numbers of the form
\begin{eqnarray} \label{eq:22}
p_j^{2\alpha_j},p_j^{2\alpha_j}\cdot p_m^{2\alpha_m},p_j^{2\alpha_j}\cdot p_m^{2\alpha_m}\cdot p_r^{2\alpha_r},... < n \\
j,m,r...=1,2,...,q, \hspace{1cm} \alpha_j\leq s_j \nonumber
\end{eqnarray}
We find the quantity of such divisors writing \eqref{eq:10} in form
\begin{equation} \label{eq:23}
n^2=\big(p_1^2\big)^{s_1}\cdot \big(p_2^2\big)^{s_2}\cdot ... \cdot \big(p_q^2\big)^{s_q}
\end{equation}
Now, considering the numbers  $p_j^2$ as indivisible and taking into account only $d_j<n$, we find, using  \eqref{eq:9}, the complete quantity of representations \eqref{eq:19} with positive and integer $k$ and $l$:
\begin{eqnarray}\label{eq:24}
N_r=\frac{(s_1+1)(s_2+1)\cdot ... \cdot (s_q+1)}{2}, \hspace{1cm} \text{if $n$ is nonsquare, i.e. some of $s_j$ are odd}
\end{eqnarray}
and
\begin{eqnarray}\label{eq:25}
N_r=\frac{(s_1+1)(s_2+1)\cdot ... \cdot (s_q+1)-1}{2}, \hspace{1cm} \text{if $n$ is square, i.e. all $s_j$ are even}
\end{eqnarray}

Hence,  the representations \eqref{eq:19} are given by expression:
\begin{equation}\label{eq:26}
n=\bigg(\frac{n+d_i}{2\sqrt{d_i}}\bigg)^2-\bigg(\frac{n-d_i}{2\sqrt{d_i}}\bigg)^2,
\end{equation}
where $di$  are the divisors of $n^2$  which are  less than $n$  and are the squares of integers (including 1).

For example, if  $n = 3465 = 3^2\cdot 5\cdot 7\cdot 11, n^2 =  3^4\cdot 5^2\cdot 7^2\cdot 11^2$, then the quantity of representations \eqref{eq:19} is equal:
\begin{equation*}
N_r=\frac{(2+1)(1+1)^3}{2}=12,
\end{equation*}
and we may calculate them on \eqref{eq:26} by $d_i$  equal to:
\begin{equation*}
1;  3^2;   5^2;   7^2 ;   3^4 ;   11^2;  3^2\cdot 5^2;   3^2\cdot 7^2;   3^2\cdot 11^2;   5^2 \cdot 7^2;   3^4\cdot 5^2; 5^2\cdot11^2
\end{equation*}

Correspondently, the 12 pairs $(k,l)$,  in  \eqref{eq:19} are: (1733, 1732);   (579, 576);   (349, 344);   (251, 244);   (197, 188);   (163,152); (123,108);   (93,72);   (69,36);   (67, 32);   (61,16);   (59,4).

\subsection{} Let $n$ be even and we want to represent it as 
\begin{equation}\label{eq:27}
n=k^2-l^2, \hspace{1cm} k,l>0
\end{equation}
According to \eqref{eq:2} written as
\begin{equation*}
2(kl)^2+n^2=(k^2+l^2)^2,
\end{equation*}
$2kl$ and $k^2-l^2$ are the elements of Pythagorean triples generated by the even number $n$. Then, according to \eqref{eq:15},
\begin{equation}\label{eq:28}
k^2+l^2=\frac{(\frac{n}{2})^2}{d_i}+d_i,
\end{equation} 
where $d_i$ are the divisors of $\big(\frac{n}{2}\big)^2$ which are less than $\frac{n}{2}$. From \eqref{eq:27} and \eqref{eq:28} we find:
\begin{equation}\label{eq:29}
k=\frac{n+2d_i}{2\sqrt{2d_i}}, \hspace{1cm} l=\frac{n-2d_i}{2\sqrt{2d_i}}
\end{equation}

It follows from \eqref{eq:29}, that $k$  and $l$ will be positive integers if $d_i<\frac{n}{2}$, and have the form 
\begin{equation} \label{eq:30}
d_i=2^{2a+1}\cdot (2b+1)^{2c}, \hspace{1cm} \text{integers} \hspace{0.5cm} a,b,c\geq 0
\end{equation}
Let’s find the quantity of such divisors, using canonic expansion \eqref{eq:13} of $\big(\frac{n}{2}\big)^2$. Since the maximal value of $2a+1$   is $2s_1-1$ , then in \eqref{eq:30} $0\leq a < s_1$, and $\big(2b+1\big)^{2c}$ are the numbers of the form \eqref{eq:22}, where $j=2,...q$.
 According to \eqref{eq:13} and  \eqref{eq:23} the quantity of all odd divisors of $\big(\frac{n}{2}\big)^2$  having the form  \eqref{eq:22} is 
 
 \begin{equation*}
 (s_2+1)\cdot (s_3+1)\cdot ... \cdot (s_q+1)
 \end{equation*}
 
Then the quantity of all even $d_i$  having the form \eqref{eq:30} is 
 
\begin{equation}\label{eq:31}
\frac{1}{2}\cdot2s_1\cdot(s_2+1)\cdot... \cdot(s_q+1)
\end{equation}

Taking into account that for obtaining of the all different representations \eqref{eq:27} only $d_i<\frac{n}{2}$ must be used, we obtain the following result.

The quantity of the representations \eqref{eq:27} of the even number $n$ is: 

\begin{align}\label{eq:32}
N_r=\frac{s_1\cdot (s_2+1)\cdot...\cdot(s_q+1)}{2}, \hspace{0.5cm}\parbox{15em}{if $n$ is nonsquare, i.e. either $s_1$ is even or some of $s_2,\cdot...\cdot,s_q$ are odd}
\end{align}         
and
\begin{equation}\label{eq:33}
N_r=\frac{s_1\cdot (s_2+1)\cdot...\cdot(s_q+1)-1}{2}, \hspace{0.5cm}\parbox{15em}{if $n$ is square, i.e. $s_1$ is odd and all $s_2,\cdot...\cdot,s_q$ are even}
\end{equation}

All these representations are given by expression
\begin{equation}\label{eq:34}
n=\bigg(\frac{n+2d_i}{2\sqrt{2d_i}}\bigg)^2-\bigg(\frac{n-2d_i}{2\sqrt{2d_i}}\bigg)^2
\end{equation}
where $d_i$  are the divisors of $\big(\frac{n}{2}\big)^2$, which are less the $\frac{n}{2}$  and have the form \eqref{eq:30}. One can see that $k$  and $l$  have the same parity.

As it follows from \eqref{eq:32}, if $s_1=0$, then $N_r=0$, i.e. the even $n$ indivisible by 4 cannot be represented as the difference of two squares. Remember that they also cannot generate the  primitive  Pythagorean triples.

In \cite{nyblom} the quantities of the representations \eqref{eq:19} and \eqref{eq:27} have been obtained by method of factorization of number $n$ in two factors. If we take the designations used there, then our expressions for $N_r$ coincide with  \cite{nyblom}, with only difference that, in contrast to \cite{nyblom}, \eqref{eq:25} and \eqref{eq:33} don’t include the trivial case  $k=\sqrt{n},l=0$.

Example.
$n= 900 = 2^2\cdot 3^2\cdot5^2;    \frac{n}{2} = 2\cdot 3^2\cdot 5^2 ; \big(\frac{n}{2}\big)^2= 2^2\cdot3^4\cdot5^4; N_r=\frac{1\cdot(2+1)^2-1}{2}=4 $. The appropriate divisors $d_i$ are $d_1=2; d_2=2\cdot3^2; d_3=2*3^4; d_4=2\cdot 5^2$, and the corresponding pairs $(k,l)$ are: (226, 224);  (78,72); (50,40); (34,16).

Now we return to the Pythagorean triples and find out how  to calculate all triples by \eqref{eq:2}.
 
Let the odd number $n=k^2-l^2$ be given.

Note that if all $s_j=1$, then the quantity $N_r$ is equal to quantity of primitive triples, i.e. $N_r=2^{q-1}$. Generally speaking $N_r\geq N_p$ , and their difference gives some  quantity of nonprimitive triples, but not all nonprimitive triples are exhausted by them. The  remaining $N_{tr}-N_r$  triples are obtained from \eqref{eq:14} by divisors $d_i<n$, which are not the squares of integer numbers. By the such $d_i$ expressions \eqref{eq:21} give the irrational $k$ and $l$, while $2kl, k^2-l^2$ and $k^2+l^2$  remain integer. On the other side, they are just those nonprimitive triples which cannot be obtained from \eqref{eq:2} with integers $k$ and $l$. Therefore the complete quantity $N_tr$ can be obtained by \eqref{eq:2}, using the  irrational $k$, $l$ given by expressions \eqref{eq:21}, alongside with integer $k$ and $l$.

Let the even number $n=2kl$ be given. Then according to \eqref{eq:15} we have 
\begin{equation*}
k^2+l^2=\frac{\big(\frac{n}{2}\big)^2}{d_i}+d_i, \hspace{1cm} k^2-l^2=\frac{\big(\frac{n}{2}\big)^2}{d_i}-d_i
\end{equation*}
and
\begin{equation}\label{eq:35}
k=\frac{n}{2\sqrt{d_i}}, \hspace{1cm} l=\sqrt{d_i}
\end{equation}

Therefore, all $N_{tr}$ triples  generated by the even number $n$,  may be obtained from \eqref{eq:2} using $k=\frac{n}{2\sqrt{d_i}}$, $l=\sqrt{d_i}$, where $d_i$ are all divisors of $\big(\frac{n}{2}\big)^2$, which are less than $\frac{n}{2}$  . 
 
\section{Conclusion}

In the more later  methods of  generation  of  the Pythagorean triples the  different special representations of generating numbers are used, in particular, the representation by Fibonacci numbers \cite{horadam1961fibonacci}, geometrical representations, e.g. Dickson's method \cite{dickson2013history}, and others \cite{amato2017characterization,2012arXiv1201.2145R}. Therefore they don’t give the universal formulas for finding all triples. In this work  we have obtained the general  formulas, giving all primitive and nonprimitive triples generated by given number which is represented in the most general form, by its canonic expansion. We have also used this method to find all representations of integers as a difference of two squares and to reveal the relation between  quantities of such representations and triples.

We have  also shown  how one can use the  Euclid’s formula \eqref{eq:2} to find all triples, which cannot be obtained by this formula with integers $k$  and $l$.

All these results are obtained with the common method using the formulas  \eqref{eq:14}, \eqref{eq:15} and canonic expansions of appropriate numbers.


\end{document}